\documentclass[reqno,11pt]{amsart}
\newlength{\myhmargin}  
\usepackage[text={6in,8in},centering,letterpaper,dvips]{geometry}
\usepackage{amsfonts,amsmath,amssymb,amsthm,verbatim,setspace,mathabx}
\usepackage{color}
\usepackage{graphicx}
\usepackage{caption}
\usepackage{subcaption}
\usepackage[none]{hyphenat}
\usepackage{hyperref}
\hypersetup{backref,colorlinks=true,citecolor=blue,linkcolor=blue}
\title{Algebraic torsion via Heegaard Floer homology}
\author[Kutluhan]{\c Ca\u gatay Kutluhan}
\author[Mati\'{c}]{Gordana Mati\'{c}}
\author[Van Horn-Morris]{Jeremy Van Horn-Morris}
\author{Andy Wand}
\address{Department of Mathematics, University at Buffalo}
\email{kutluhan@buffalo.edu}
\address{Department of Mathematics, University of Georgia}
\email{gordana@math.uga.edu}
\address{Department of Mathematics, University of Arkansas}
\email{jvhm@uark.edu}
\address{School of Mathematics and Statistics, University of Glasgow}
\email{andy.wand@glasgow.ac.uk}
\thanks{\c Ca\u gatay Kutluhan was supported in part by NSF grant DMS-1360293. \\ Gordana Mati\'c and Jeremy Van Horn-Morris were supported in part by Simons Foundation grants No.~246461 and No.~279342, respectively. \\Andy Wand was supported in part by ERC grant \textsc{geodycon}.}

\theoremstyle{definition}
\newtheorem*{defn}{Definition}
\newtheorem*{rmk}{Remark}

\newtheorem*{ack}{Acknowledgments}
\begin{document}
\bibliographystyle{amsalpha}
\maketitle
\begin{abstract}
This is a prelude to the authors' work in \cite{KMVHMW2} in which they define an invariant of contact structures in dimension three as a refinement of the contact invariant in Heegaard Floer homology \cite{OzsvathSzabo4}. The idea is to port Hutchings's recipe that produces an ECH analog of Latschev and Wendl's algebraic $k$-torsion (see \cite[Appendix]{LatschevWendl}) to Heegaard Floer homology using the construction of an isomorphism between the latter and Seiberg--Witten Floer homology by Lee, Taubes, and the first  author (see \cite{KLT} for an outline). This isomorphism factors through a variant of Hutchings's embedded contact homology, denoted \emph{ech}. This note outlines Hutchings's recipe in the context of \emph{ech} and explains how it translates into Heegaard Floer homology.
\end{abstract}
\section*{An analog of algebraic $k$-torsion}
To set the stage, let $M$ be a closed, connected, and oriented $3$-manifold endowed with a co-oriented contact structure $\xi$. It is understood that the orientation on $M$ is induced by $\xi$. Fix an abstract open book decomposition $(S,\phi)$ of $M$ supporting $\xi$. Here, $S$ is a compact oriented surface of genus $g$ with $\textsc{b}$ boundary components, called the \emph{page}, and $\phi$ is an orientation preserving diffeomorphism of $S$ which restricts to identity in a neighborhood of the boundary, called the \emph{monodromy}. The manifold $M$ is diffeomorphic to $S\times[0,1]/\sim$ where $(p,1)\sim(\phi(p),0)$ for any $p\in S$ and $(p,t)\sim(p,t')$ for any $p\in\partial S$ and $t,t'\in[0,1]$. Let $\textsc{g}=2g+\textsc{b}-1$, and fix a self-indexing Morse function on $S$ that has a single maximum and attains its minimum along $\partial S$. Then a suitable pseudo-gradient vector field for the latter defines a basis of arcs $\{a_1,\dots,a_\textsc{g}\}$ on $S$, that is, a pairwise disjoint collection of properly embedded arcs cutting $S$ into a polygon. This basis together with the monodromy $\phi$ defines a Heegaard diagram $(\Sigma,\{\beta_1,\dots,\beta_\textsc{g}\},\{\alpha_1,\dots,\alpha_\textsc{g}\})$ for $-M$ as in \cite[\S 3.1]{HKM}. To be more explicit, let $\{b_1,\dots,b_\textsc{g}\}$ be a collection of arcs on $S$ where $b_i$ is isotopic to $a_i$ via a small isotopy satisfying the following conditions:
\begin{itemize}\leftskip-0.25in
\item The endpoints of $b_i$ are obtained from the endpoints of $a_i$ by pushing along $\partial S$ in the direction of the boundary orientation,
\item $a_i$ intersects $b_i$ transversally at one point in the interior of $S$,
\item Having fixed an orientation of $a_i$, there is an induced orientation on $b_i$, and the sign of the oriented intersection $a_i\cap b_i$ is positive.
\end{itemize}
Then $\Sigma= S\times\{\frac{1}{2}\}\cup_{\partial S}-S\times\{0\}$, $\alpha_i=a_i\times\{\frac{1}{2}\}\cup a_i\times\{0\}$, and $\beta_i=b_i\times\{\frac{1}{2}\}\cup\phi(b_i)\times\{0\}$. Alternatively, the pseudo-gradient vector field on $S$ and the monodromy $\phi$ can be used to define a self-indexing Morse function $f$ and a pseudo-gradient vector field $v$ on $-M$ which yield the desired Heegaard diagram. Note that there is a natural pairing of the index-1 and index-2 critical points indicated by the labeling of the $\alpha$- and $\beta$-curves. Now fix a basepoint $z\in\Sigma\smallsetminus\bigcup_{i\in\{1,\dots,\textsc{g}\}}(\alpha_i\cup\beta_i)$ according to the convention in \cite[\S 3.1]{HKM}. By changing the monodromy in its isotopy class, one can make sure that the pointed Heegaard diagram $(\Sigma,\{\beta_1,\dots,\beta_\textsc{g}\},\{\alpha_1,\dots,\alpha_\textsc{g}\},z)$ is strongly admissible for any given $Spin^c$ structure. Then, as in \cite[Lemma 1.1]{KLT2}, there exists an area form $w_\Sigma$ on $\Sigma$ with total area equal to $2$, and the signed area of each periodic domain $\mathcal{P}$ on this Heegaard diagram is equal to the pairing $\langle c_1(\mathfrak{s}_\xi), \mathcal{H}(\mathcal{P})\rangle$, where $\mathcal{H}(\mathcal{P})\in H_2(M;\mathbb{Z})$ is the homology class corresponding to $\mathcal{P}$. 

Next, following \cite[\S 1a]{KLT2}, construct a new manifold $Y$ out of $M$ by adding $\textsc{g}+1$ $1$-handles such that one, denoted $\mathcal{H}_0$, is attached along the index-0 and index-3 critical points of $f$, whereas the remaining $1$-handles are attached along pairs of index-1 and index-2 critical points of $f$ as prescribed by the pairing of the $\alpha$ and $\beta$ curves. The latter kind of handles are denoted by $\mathcal{H}_{i}$ for $i\in\{1,\dots,\textsc{g}\}$. The resulting manifold $Y$ is diffeomorphic to $M\#_{\textsc{g}+1}S^1\times S^2$ via an orientation preserving diffeomorphism. With the preceding understood, use $v$ and $w_\Sigma$, as in \cite[\S 1a-\S 1d]{KLT2}, to construct a stable Hamiltonian structure on $Y$. The latter is a pair $(a,w)$ where $a$ is a smooth $1$-form and $w$ is a smooth closed $2$-form such that $da=hw$ for some smooth function $h:Y\to\mathbb{R}$, and $a\wedge w$ is nowhere zero. Associated to this stable Hamiltonian structure is a vector field $R$ satisfying $w(R,\cdot)=0$ and $a(R)=1$, which agrees with the pseudo-gradient vector field $v$ on $Y\smallsetminus \bigcup_{i\in\{0,1,\dots,\textsc{g}\}}\mathcal{H}_{i}$. Of interest are periodic orbits of the vector field $R$, called the \emph{Reeb vector field}. A complete description of periodic orbits of the Reeb vector field is provided in \cite[\S 2]{KLT2}. In particular, there is a unique embedded periodic orbit passing through the basepoint $z$, denoted by $\gamma_z$. This orbit intersects every cross-sectional sphere of $\mathcal{H}_0$ exactly once. Meanwhile, for each $i\in\{1,\dots,\textsc{g}\}$, there are exactly two distinct embedded periodic orbits inside $\mathcal{H}_i$, denoted by $\gamma_{i}^+$ and $\gamma_{i}^-$, which are hyperbolic and homologically trivial. 

Periodic orbits of the Reeb vector field and pseudo-holomorphic curves in $\mathbb{R}\times Y$ for suitably generic almost complex structure are used in \cite[Appendix A]{KLT2} to define a variant of Hutchings's embedded contact homology for certain homology classes $\Gamma\in H_1(Y;\mathbb{Z})$. First, we recall some basic definitions. An \emph{orbit set} is a finite collection of the form $\{(\gamma,m)\}$ where $\gamma$ are distinct embedded periodic orbits of the Reeb vector field and $m$ is a positive integer. An orbit set $\{(\gamma,m)\}$ with the extra requirement that $m=1$ when $\gamma$ is hyperbolic is called an \emph{admissible orbit set}. The homology class of an orbit set $\Theta=\{(\gamma,m)\}$ is defined by $\sum m[\gamma]$. Of interest here are admissible orbit sets with homology class $\Gamma\in H_1(Y;\mathbb{Z})$ such that $\langle PD(\Gamma),[S^2]\rangle=0$ if $[S^2]$ is a generator of $H_2(\mathcal{H}_0;\mathbb{Z})$, $\langle PD(\Gamma),[S^2]\rangle=1$ if $[S^2]$ is the positive generator of $H_2(\mathcal{H}_i;\mathbb{Z})$ for $i\in\{1,\dots,\textsc{g}\}$, and $\mathfrak{s}-\mathfrak{s}_\xi=PD(\Gamma|_M)$ where $\mathfrak{s}_\xi$ is the canonical $\mathrm{Spin}^c$ structure on $M$ defined by the contact structure $\xi$. Note that the homology class $\Gamma$ is uniquely determined by these properties once a $\mathrm{Spin}^c$ structure $\mathfrak{s}$ on $M$ is fixed. For such $\Gamma\in H_1(Y;\mathbb{Z})$, it follows from \cite[\S 2]{KLT2} that any $\gamma$ that belongs to an admissible orbit set with homology class $\Gamma$ is hyperbolic. Hence, we may drop the integer $m$ from the notation. The free $\mathbb{Z}$-module generated by admissible orbit sets with homology $\Gamma$ will be denoted by $\widehat{\mathit{ecc}}(Y,\Gamma)$. With the preceding understood, the relationship between admissible orbit sets and Heegaard Floer generators is described in \cite[Proposition 2.8]{KLT2}. To be more precise, the set of admissible orbit sets associated to a Heegaard Floer generator is in 1--1 correspondence with the set $\prod_{i\in\{1,\dots,\textsc{g}\}}(\mathbb{Z}\times\textsc{o})$ where $\textsc{o}=\{0,+1,-1,\{+1,-1\}\}$.

Having fixed a generic almost complex structure $J$ on $\mathbb{R}\times Y$ satisfying the conditions in \cite[\S 3a]{KLT2}, the differential $\widehat{\partial}_{\mathit{ech}}$ on $\widehat{\mathit{ecc}}(Y,\Gamma)$ is defined to be the endomorphism of $\widehat{\mathit{ecc}}(Y,\Gamma)$ sending a generator $\Theta_+$ to
\[\sum \sigma(\Theta_+,\Theta_-)\Theta_-,\]
where $\sigma(\Theta_+,\Theta_-)$ is a signed count, modulo $\mathbb{R}$-translation, of ECH index-$1$ $J$-holomorphic curves asymptotic to $\Theta_+$ at $+\infty$ and to $\Theta_-$ at $-\infty$, and disjoint from $\mathbb{R}\times\gamma_z$. Note that $\mathbb{R}\times\gamma_z$ is a $J$-holomorphic cylinder. Hence, by positivity of intersections of pseudo-holomorphic curves, $\mathbb{R}\times\gamma_z$ introduces a filtration on the chain complex defined in \cite[\S 1b]{KLT3}, much like the filtration dictated by the basepoint in Heegaard Floer homology. It follows from \cite{KLT3} (cf. \cite{HutchingsTaubes, HutchingsTaubes2}) that $\widehat{\partial}_{\mathit{ech}}\circ\widehat{\partial}_{\mathit{ech}}=0$, and
\[\widehat{\mathit{ech}}(Y,\Gamma):=H_\ast(\widehat{\mathit{ecc}}(Y,\Gamma),\widehat{\partial}_{\mathit{ech}})\cong\widehat{\mathit{HF}}(-M,\mathfrak{s}_\xi+PD(\Gamma|_M))\otimes\mathbb{Z}^{2^\textsc{g}}.\]
A detailed description of $J$-holomorphic curves in $\mathbb{R}\times Y$ is given in \cite[\S 3 and \S 4]{KLT2}.

Given orbit sets $\Theta_+$ and $\Theta_-$ with the same homology class, Hutchings defines variants of his ECH index for relative homology classes in $H_2(Y,\Theta_+,\Theta_-)$, denoted by $J_0$, $J_+$, and $J_-$ (see \cite[\S 6]{Hutchings3}). Among other things, these are additive in the sense that for $Z_1\in H_2(Y,\Theta_+,\Theta)$ and $Z_2\in H_2(Y,\Theta,\Theta_-)$, we have
\[J_\circ(Z_1+Z_2)=J_\circ(Z_1)+J_\circ(Z_2),\]
where $\circ\in\{0,+,-\}$ (see \cite[Proposition 6.5]{Hutchings3}). If $C$ is an embedded $J$-holomorphic curve in $\mathbb{R}\times Y$ with positive ends at an admissible orbit set $\Theta_+$ and negative ends at an admissible orbit set $\Theta_-$, then $C$ has ends at distinct embedded periodic orbits of the Reeb vector field, and hence
 \[J_+(C)=-\chi(C)+|\Theta_+|-|\Theta_-|,\]
where $|\cdot|$ denotes the cardinality of an admissible orbit set. This is due to the fact that all relevant periodic orbits of the Reeb vector field are hyperbolic. We can rewrite the above formula as
\begin{equation}\label{eq:jplus}
J_+(C)=\sum_{C_j\subset C} (2g_j-2 +2|\Theta_{j+}|),
\end{equation}
where each $C_j$ denotes a connected component of $C$, $g_j$ denotes the genus of $C_j$, and each $\Theta_{j+}\subset\Theta_+$ is the collection of periodic orbits of the Reeb vector field that comprise the positive ends of $C_i$. Clearly, $2\,|\,J_+(C)$. Note that no subcollection of periodic orbits from an admissible orbit set $\Theta$ can be homologically trivial unless $\Theta$ contains one or more of the homologically trivial $\gamma_{i}^+$, $\gamma_{i}^-\subset\mathcal{H}_{i}$, and admissible orbit sets whose restriction to each $\mathcal{H}_i$ correspond to the pair $(0,0)\in\mathbb{Z}\times\textsc{o}$ generate a subcomplex, which we denote by $\widehat{\mathit{ecc}}_\circ(Y,\Gamma)$. If $\Theta_+$ and $\Theta_-$ are generators of $\widehat{\mathit{ecc}}_\circ(Y,\Gamma)$, then every connected component of an embedded $J$-holomorphic curve $C$ with positive ends at $\Theta_+$ and negative ends at $\Theta_-$ has non-empty positive and negative ends, and \eqref{eq:jplus} implies that $J_+(C)\geq 0$. Following \cite[Appendix]{LatschevWendl}, we decompose the differential $\widehat{\partial}_{ech}$ as
\[\widehat{\partial}_{ech}=\partial_0+\partial_1+\cdots+\partial_k+\cdots,\]
where $\partial_k$ counts admissible ECH index-1 $J$-holomorphic curves $C$ with $J_+(C)=2k$ and empty intersection with $\mathbb{R}\times\gamma_z$. Since $J_+$ is additive under gluing of $J$-holomorphic curves (see \cite[Proposition 6.5a]{Hutchings3}), the above decomposition induces a spectral sequence with pages $E^k(S,\phi,a;J)$.

With the preceding understood, let $\mathfrak{s}=\mathfrak{s}_\xi$ be the canonical $Spin^c$ structure on $M$. The set of points $\mathbf{x}_\xi=\{x_1,\dots,x_\textsc{g}\}$ that defines the Ozsv\'ath--Szab\'o contact class in the Honda--Kazez--Mati\'{c} description \cite{HKM} specifies a distinguished generator $\Theta_\xi$ of $\widehat{\mathit{ecc}}_\circ(Y,\Gamma)$ which is a disjoint union of $\textsc{g}$ embedded periodic orbits of the Reeb vector field each representing a positive generator of $H_1(S^1\times S^2;\mathbb{Z})$. 

\begin{defn}
Let $AT(S,\phi,a;J)$ be the smallest non-negative integer $k$ such that the orbit set $\Theta_\xi$ represents the trivial class in $E^{k+1}(S,\phi,a;J)$. Borrowing the terminology from \cite{LatschevWendl}, we say that \emph{$(M,\xi)$ has algebraic $k$-torsion} if $AT(S,\phi,a;J)=k$ for some choice of $(S,\phi,a)$ and generic almost complex structure $J$.
\end{defn}

We conclude our discussion by giving a description of the relative filtration by $J_+$, and hence, a description of algebraic $k$-torsion, in terms of the data on the Heegaard diagram associated to $(S,\phi,a)$. In order to do so, we use the isomorphism between Heegaard Floer homology and \emph{ech} \cite{KLT2, KLT3}. The construction of this isomorphism exploits the cylindrical reformulation of Heegaard Floer homology due to Lipshitz \cite{Lipshitz}. In this regard, given an almost complex structure $J_{\textit{HF}}$ on $\Sigma\times [0,1]\times\mathbb{R}$ satisfying conditions (\textbf{J1})--(\textbf{J5}) in \cite[\S 1]{Lipshitz}, there exists an almost complex structure $J$ on $\mathbb{R}\times Y$ satisfying the conditions in \cite[\S 3a]{KLT2}. In the present context, where we restrict ourselves to the subcomplex $\widehat{\mathit{ecc}}_\circ(Y,\Gamma)$, there is a 1--1 correspondence between ECH index-$1$ $J$-holomorphic curves in $\mathbb{R}\times Y$ with ends at generators of this subcomplex and Fredholm index-$1$ $J_{\mathit{HF}}$-holomorphic curves in $\Sigma\times [0,1]\times\mathbb{R}$ that satisfy conditions~(\textbf{M0})--(\textbf{M6}) in \cite[\textsection 1]{Lipshitz}. Therefore, $(\widehat{\mathit{ecc}}_\circ(Y,\Gamma),\widehat{\partial}_{ech})$ and $(\widehat{\mathit{CF}}(-M,\mathfrak{s}),\widehat{\partial}_{\mathit{HF}})$ are canonically isomorphic as chain complexes via the aforementioned 1--1 correspondence between generators  (see \cite{KLT3} for details). 

Now let $C_L$ be a Fredholm index-$1$ $J_{\mathit{HF}}$-holomorphic curve in $\Sigma\times [0,1]\times\mathbb{R}$ satisfying conditions (\textbf{M0})--(\textbf{M6}) in \cite[\textsection 1]{Lipshitz}, and $C\subset\mathbb{R}\times Y$ be the corresponding ECH index-$1$ $J$-holomorphic curve with positive ends at an admissible orbit set $\Theta_+$ and negative ends at an admissible orbit set $\Theta_-$, both of which are generators of $\widehat{\mathit{ecc}}_\circ(Y,\Gamma)$. Due to the aforementioned 1--1 correspondence, we can write $J_+(C_L)$ for $J_+(C)$. Topologically, $C$ is obtained from $C_L$ by adding 2-dimensional $1$-handles in such a way that for each $i=1,\dots,\textsc{g}$ a $1$-handle is attached along a point in $\alpha_i$ and a point in $\beta_i$. Therefore, $-\chi(C)=-\chi(C_L)+\textsc{g}$, and 
\begin{equation}
\label{eq:j-plus-lip}
J_+(C_L)=-\chi(C_L)+\textsc{g}+|\Theta_+|-|\Theta_-|.
\end{equation} 
Lipshitz proves that $\chi(C_L)$ is determined by the relative homology class $[C_L]\in \widehat{\pi}_2(\mathbf{x}_+,\mathbf{x}_-)$ (see \cite[Proposition 4.2]{Lipshitz}, cf. \cite[Proposition $4.2'$]{Lipshitz2}); more specifically,
\begin{equation}
\label{eq:euler-lip}
\chi(C_L)=\textsc{g}-n_{\mathbf{x}_+}(\mathcal{D}([C_L]))-n_{\mathbf{x}_-}(\mathcal{D}([C_L]))+e(\mathcal{D}([C_L])).
\end{equation} 
Here, $\mathcal{D}([C_L])$ denotes the domain in the Heegaard diagram associated to $C_L$, $\mathbf{x}_+$ denotes the Heegaard Floer generator that corresponds to $\Theta_+$, $\mathbf{x}_-$ denotes the Heegaard Floer generator that corresponds to $\Theta_-$, $n_p(\mathcal{D}([C_L]))$ denotes the \emph{point measure}, namely, the average of the coefficients of $\mathcal{D}([C_L])$ for the four regions with corners at $p\in\alpha_i\cap\beta_j$, and $e(\mathcal{D}([C_L]))$ is the \emph{Euler measure} of $\mathcal{D}([C_L])$. Meanwhile, keeping the labeling of the $\alpha$ and $\beta$ curves in mind, each Heegaard Floer generator yields an element of the symmetric group, and $|\Theta_\pm|$ is equal to the number of cycles in the element of that symmetric group associated to $\mathbf{x}_\pm$, respectively. We will denote the latter quantity also by $|\cdot|$. With the preceding understood, combine \eqref{eq:j-plus-lip} and \eqref{eq:euler-lip} to obtain the following formula:
\begin{equation}\label{eq:jplus-12}
J_+([C_L])=n_{\mathbf{x}_+}(\mathcal{D}([C_L]))+n_{\mathbf{x}_-}(\mathcal{D}([C_L]))-e(\mathcal{D}([C_L]))+|\mathbf{x}_+|-|\mathbf{x}_-|.
\end{equation}
Note that the right hand side of \eqref{eq:jplus-12} depends only on $\mathcal{D}([C_L])$, and in particular, since a class in $\widehat{\pi}_2(\mathbf{x}_+,\mathbf{x}_-)$ is specified by its domain, $J_+([C_L])$ depends only on the relative homology class $[C_L]$. In fact, we may define $J_+$ for any domain $\mathcal{D}$ in the Heegaard diagram by the right hand side of \eqref{eq:jplus-12}, which by \cite[Corollary 4.10]{Lipshitz} (cf.  \cite[Proposition $4.8'$]{Lipshitz2}) could also be written as 
\begin{equation}\label{eq:jplus-prod}
J_+(\mathcal{D}):=\mu(\mathcal{D})-2e(\mathcal{D})+|\mathbf{x}_+|-|\mathbf{x}_-|,
\end{equation}
where $\mu(\mathcal{D})$ is the \emph{Maslov index} of $\mathcal{D}$. 

\begin{ack} This project was initiated at the ``\emph{Interactions between contact symplectic topology and gauge theory in dimensions 3 and 4}" workshop at Banff International Research Station (BIRS) in 2011. We would like to thank BIRS and the organizers of that workshop for their hospitality and for creating a wonderful environment for collaboration. We also thank Michael Hutchings for helpful correspondence.
\end{ack}

\begin{rmk}
We were informed that John Baldwin and David Shea Vela-Vick have independently been working on a refinement of the Heegaard Floer contact invariant, similar in spirit to what is outlined in this note. 
\end{rmk}

\bibliography{biblio}
\end{document}